\documentclass{article}
\usepackage{amsfonts, cite}
\usepackage{amsmath,amssymb,amsthm,epic,eepic,eepicemu,graphicx}
\usepackage{graphicx}
\usepackage{epstopdf}

\def\a{\alpha}
\def\b{\beta}

\newcommand{\beq}{\begin{equation}}
\newcommand{\eeq}{\end{equation}}
\newcommand{\dpl}{\displaystyle}
\newcommand{\sech}{\mathop{\rm sech}\nolimits}

\newcommand{\cR}{{\cal R}}

\def\I{{\rm I}}
\def\II{{\rm II}}
\def\III{{\rm III}}
\def\IV{{\rm IV}}
\def\V{{\rm V}}

\def\({\left(}
\def\){\right)}
\def\<{\langle}
\def\>{\rangle}

\begin{document}

\title{Yang-Baxter maps: dynamical point of view}

\date{}

\maketitle

\begin{center}

{\bf A.P.Veselov }

\bigskip

{\it Department of Mathematical Sciences, Loughborough University,
Loughborough, Leicestershire, LE 11 3TU, UK
}

\bigskip

{\it Landau Institute for Theoretical Physics, Kosygina 2,

 Moscow, 117940,  Russia

\bigskip

e-mail: A.P.Veselov@lboro.ac.uk,
}

\end{center}

\bigskip

\bigskip

{\small  {\bf Abstract}
A review of some recent results on the dynamical theory of the Yang-Baxter maps (also known as set-theoretical solutions to the quantum Yang-Baxter equation) is given.
The central question is the integrability of the transfer dynamics.
The relations with matrix factorisations, matrix KdV solitons, Poisson Lie groups, geometric crystals and tropical combinatorics are discussed and demonstrated on several concrete examples.}

\bigskip

\section*{Introduction}

The story of Yang-Baxter maps is usually traced back to Drinfeld \cite{D}, who
suggested to study the set-theoretical solutions to the quantum Yang-Baxter equation, but the first interesting examples were constructed already by Sklyanin \cite{Skl}. Drinfeld's question was immediately appreciated by Weinstein and Xu \cite{WX} who developed an approach to this problem based on the theory of Poisson Lie groups and symplectic groupoids.
Etingof, Schedler and Soloviev  \cite{ESS, Solo}, partly inspired by Hietarinta's work \cite{Hi}, took more algebraic point of view (see also \cite{ Buch, ESG, LYZ, GB, G-I, CS}).

Very interesting examples of such solutions have appeared in the theory of geometric crystals (Berenstein-Kazhdan \cite{BK} and Etingof \cite{E}). This opened remarkable links with the theory of soliton cellular automata ("box-ball" systems) and tropical combinatorics extensively developed by the Japanese school \cite{TS,  HKT, HHIKTT, HKOT, KNY, Kir, Y, KOTY, NouY}.

The study of the dynamical aspects of this problem was initiated in the paper \cite{V1},
where the term "Yang-Baxter map" was introduced. My motivation came from the soliton theory,
in particular from the theory of matrix KdV equation \cite{GV}. I was puzzled by what kind of information one can learn from the interaction of two matrix solitons and what is the role of the Yang-Baxter relation here. Analogy with the theory of solvable models in statistical mechanics and quantum inverse scattering method \cite{Yang, B, Gaudin, TF, JM}, where the quantum Yang-Baxter equation has its origin,  naturally led to the notion of the transfer dynamics and to the relations with the theory of integrable systems.  The main idea of \cite{V1} was that the investigation of the transfer dynamics should reveal the structures encapsulated in a given Yang-Baxter map. A very good example is the Adler map (see below) from the theory of the dressing chain \cite{VS, A}. The corresponding transfer dynamics is nothing else but the discrete version of the finite-gap KdV dynamics (Novikov equations \cite{Nov, DMN}) with its rich algebro-geometric and symplectic structures. 

In this extended version of the lectures given at Kyoto meeting "Combinatorial
Aspects of Integrable Systems" in July 2004 I am going to review the dynamical theory of the Yang-Baxter maps and illustrate it on several concrete examples following mainly to \cite{V1, GV, SV, RV}.

We start with the general discussion of the Yang-Baxter maps and
related transfer dynamics. The main property of the transfer maps, 
which they share with the transfer matrices, is the
pairwise commutativity. If the maps are polynomial or rational the
commutativity relation is quite strong and sometimes implies the
solvability of the corresponding dynamical system (see \cite{V}
for the discussion of some known results in this direction).
Although it is not clear what is the transfer dynamics in general,
our main motive is that it should be "integrable" 
and one can think of it also as a set-theoretical analogue of the finite-gap theory.

We illustrate this by two constructions of the Yang-Baxter equation coming from the theory of integrable systems. The first one is related to matrix factorisations and QR-type procedure, which is a discrete analogue of the Lax representation.  We explain how the Lax matrices can be read off
the map itself following the procedure suggested in \cite{SV}.

The second class is given by the interaction of solitons with non-trivial internal parameters (e.g. matrix KdV solitons). It is a well-known phenomenon in soliton theory (see e.g. \cite{Kul,Sol,G,Tsu, APT}) that the interaction of $n$ solitons is completely determined by their pairwise
interactions (so there are no multiparticle effects). The fact that the final result is independent of the order
of collisions means that the map determining the interaction of two solitons satisfies the Yang-Baxter relation. 

We briefly consider also the examples of the Yang-Baxter maps related to geometric crystals
and tropical combinatorics \cite{E, KNY, KOTY}. Then we discuss the Hamiltonian theory of the Yang-Baxter maps following recent work by Reshetikhin and the author \cite{RV} partly based on some general considerations from \cite{WX}. The main observation is that if the Lax matrices form symplectic leaves in some Poisson Lie groups, then the corresponding transfer dynamics is symplectic. 
In the last section we discuss the classification problem of the Yang-Baxter maps and its relation with the classification of the quadrirational maps  done recently by Adler, Bobenko and Suris \cite{ABS2}.

\section*{Yang-Baxter maps and transfer dynamics}

The original {\it quantum Yang-Baxter equation} is the following relation on a linear operator $R : V \otimes V \rightarrow V \otimes V:$
\begin{equation}
\label{YB}
R_{12} R_{13} R_{23} = R_{23} R_{13} R_{12},
\end{equation}
where $R_{ij}$ is acting in $i$-th and $j$-th components of the tensor product $V \otimes V \otimes V$
(see e.g. \cite{JM}).

Following Drinfeld's suggestion \cite{D} consider the following set-theoretical version of this equation.

Let $X$ be any set, $R: X \times X \rightarrow X \times X$ be a map from its square into itself.
Let $R_{ij}: X^{n} \rightarrow X^{n}, \quad X^{n} = X \times X \times .....\times X$
be the maps which acts as $R$ on $i$-th and $j$-th factors and identically on the others.
More precisely, if $R(x, y) = (f(x,y), g(x,y)), x,y \in X$ then
$R_{ij} (x_1, x_2, \dots, x_n) = (x_1, x_2, \dots, x_{i-1}, f(x_i,x_j), x_{i+1}, \dots, x_{j-1}, g(x_i,x_j), x_{j+1},
\dots,x_n)$
for $i<j$ and
$(x_1, x_2, \dots, x_{j-1}, g(x_i,x_j), x_{j+1},\dots, x_{i-1}, f(x_i,x_j), x_{i+1},\dots,x_n)$
otherwise.
In particular for $n=2$
$R_{21}(x,y) = (g(y,x), f(y,x)).$ If $P: X^2 \rightarrow X^2$ is the permutation of $x$ and $y$: $P(x,y) = (y,x)$, then
obviously we have $$R_{21} = P R P.$$

Following \cite{V1} we will call $R$ {\it Yang-Baxter map}
if it satisfies the Yang-Baxter relation (\ref{YB}) considered as
the equality of the maps of $X \times X \times X$ into itself. If
additionally $R$ satisfies the relation
\begin{equation}
\label{U}
R_{21} R = Id,
\end{equation}
we will call it {\it reversible Yang-Baxter map}.

The condition (\ref{U}) is usually called the {\it unitarity} condition, but in our case
the term {\it reversibility} is more appropriate since in dynamical systems terminology
this condition means that the map $R$
is reversible with respect to the permutation $P$.

Note that if we consider the linear space $V = {\bf C}^{X}$
spanned by the set $X$, then any Yang-Baxter map $R$ induces a
linear operator in $V \otimes V$ which satisfies the quantum
Yang-Baxter equation in the usual sense. Therefore we indeed have
a very special class of solutions to this equation. If $X$ is a
finite set and $R$ is a bijection then we have permutation-type
solutions discussed in \cite{Hi}. However this point of view on
the Yang-Baxter maps seems to be artificial because it does not
reflect the nature of the maps (like birationality).

{\bf Remark.} We should warn the reader that some authors (see
e.g. \cite{Buch, E, CS}) consider instead of the map $R$ a closely
related map $S = P R,$ where $P$ is the permutation map. The
corresponding maps satisfies instead of the Yang-Baxter relation
the {\it braid relation}
\begin{equation}
\label{braid} S_{12}S_{23}S_{12} = S_{23}S_{12}S_{23}
\end{equation}
and instead of reversibility the {\it involutivity} condition
\begin{equation}
\label{inv} S^2 = Id
\end{equation}
One should be careful here since the groups generated by the
$R_{ij}$ and $S_{ij}$ maps could be different. For example, the
group generated by $R$ in general is ${\bf Z}$ while for $S$ it is ${\bf
Z}_2.$ This is an important difference from dynamical point of
view (see below).


We can represent the relations (\ref{YB}, \ref{U}) in the standard diagrammatic way as shown on the Fig. 1 and 2.

\begin{figure}[htbp]
\begin{center}
\includegraphics[width=10cm]{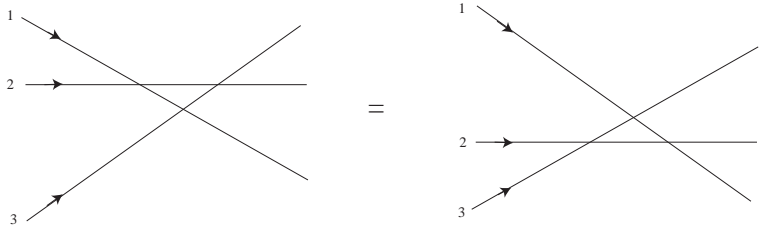}
\caption{Yang-Baxter relation}\label{fig:YBt}
\end{center}
\end{figure}


\begin{figure}[htbp]
\begin{center}
\includegraphics[width=10cm]{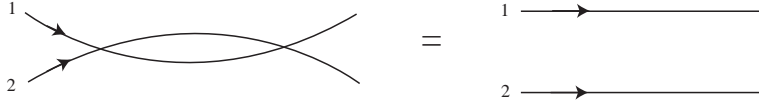}
\caption{Reversibility}\label{fig:Rev}
\end{center}
\end{figure}

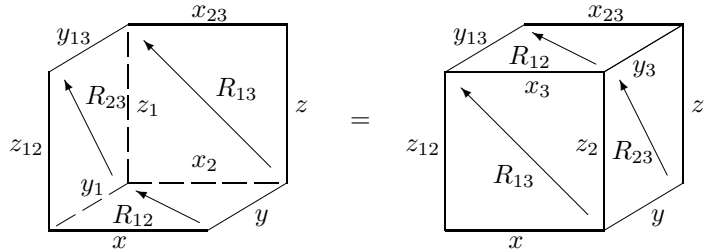
\begin{figure}[bp]
\setlength{\unitlength}{0.06em}
\begin{center}
\begin{picture}(450,170)(-30,-20)
  \put( 0,  0){\line(1,0){100}}
  \put(100,  0){\line(5,3){50}}
  \put(150,30){\line(0,1){100}}
  \put(50,130){\line(1,0){100}}
  \multiput(50,30)(20,0){5}{\line(1,0){15}}
  \put(  0, 0){\line(0,1){100}}
  \multiput(50,30)(0,20){5}{\line(0,1){15}}
  \put(  0,100){\line(5,3){50}}
  \multiput(50,30)(-16.67,-10){3}{\line(-5,-3){12}}
  \put(40,35){\vector(-1,2){30}}
  \put(140,40){\vector(-1,1){80}}
  \put(95,5){\vector(-2,1){40}}
  \put(105,85){$R_{13}$}
  \put(23,80){$R_{23}$}
  \put(40,5){$R_{12}$}
  \put(40,-11){$x$}
  \put(130,6){$y$}
  \put(155,75){$z$}
  \put(90,135){$x_{23}$}
  \put(5,120){$y_{13}$}
  \put(-25,50){$z_{12}$}
  \put(90,40){$x_2$}
  \put(20,25){$y_1$}
  \put(55,75){$z_1$}
  \put(190,65){$=$}
  \put(250,  0){\line(1,0){100}}
  \put(250, 0){\line(0,1){100}}
  \put(250,100){\line(1,0){100}}
  \put(350, 0){\line(0,1){100}}
  \put(350, 0){\line(5,3){50}}
  \put(400,30){\line(0,1){100}}
  \put(350,100){\line(5,3){50}}
  \put(300,130){\line(1,0){100}}
  \put(250,100){\line(5,3){50}}
  \put(350,100){\line(5,3){50}}
  \put(390,35){\vector(-1,2){30}}
  \put(340,10){\vector(-1,1){80}}
  \put(345,105){\vector(-2,1){40}}
  \put(280,30){$R_{13}$}
  \put(355,45){$R_{23}$}
  \put(290,105){$R_{12}$}
  \put(290,-11){$x$}
  \put(380,6){$y$}
  \put(405,75){$z$}
  \put(340,135){$x_{23}$}
  \put(255,120){$y_{13}$}
  \put(225,50){$z_{12}$}
  \put(300,87){$x_3$}
  \put(368,100){$y_3$}
  \put(332,50){$z_2$}
\end{picture}
\end{center}
\caption{``Cubic'' representation of the Yang--Baxter relation}
\label{Fig:YB}
\end{figure}

One should mention that there exist another (dual) way to
represent the Yang-Baxter relation,  which shows the
relation with the so-called {\it 3D consistency condition} in the
theory of integrable discrete equations on quad-graphs
\cite{BS, N, ABS}. In this representation the fields (elements of $X$)
are assigned to the edges of elementary quadrilaterals rather than
to the vertices. The corresponding "cubic" representation of the
Yang-Baxter relation is illustrated in the Fig. \ref{Fig:YB}.

The left--hand side of (\ref{YB})
corresponds to the chain of maps along the three rear faces of the
cube on Fig. \ref{Fig:YB}:
\[
R_{12}:(x,y)\mapsto(x_2,y_1),\quad
R_{13}:(x_2,z)\mapsto(x_{23},z_1),\quad
R_{23}:(y_1,z_1)\mapsto(y_{13},z_{12}),
\]
while its right--hand side corresponds to the chain of the maps
along the three front faces of the cube:
\[
R_{23}:(y,z)\mapsto(y_3,z_2),\quad
R_{13}:(x,z_2)\mapsto(x_3,z_{12}),\quad R_{12}:(x_3,y_3)\mapsto
(x_{23},y_{13}).
\]
So, (\ref{YB}) assures that two ways of obtaining
$(x_{23},y_{13},z_{12})$ from the initial data $(x,y,z)$ lead to
the same results.

One can introduce also a more general parameter-dependent Yang-Baxter
equation as the relation
\begin{equation}
\label{sYB}
R_{12}(\lambda_1, \lambda_2) R_{13}(\lambda_1, \lambda_3) R_{23} (\lambda_2, \lambda_3) = R_{23}(\lambda_2, \lambda_3)
R_{13}(\lambda_1, \lambda_3) R_{12}(\lambda_1, \lambda_2)
\end{equation}
and the corresponding unitarity (reversibility) condition as
\begin{equation}
\label{sU}
R_{21}(\mu,\lambda) R(\lambda, \mu) = Id.
\end{equation}
Although it can be reduced to the usual case by considering
$\tilde X = X \times {\bf C}$ and $\tilde R (x,\lambda; y, \mu) =
R (\lambda,\mu) (x,y)$ often it is useful to keep the parameter
separately (see the examples below).

In the theory of the quantum Yang-Baxter equation the following {\it transfer matrices}
$t^{(n)}: V^{\otimes n} \rightarrow V^{\otimes n}$ play a fundamental role.
They are defined as the trace of the monodromy matrix
\begin{equation}
\label{transfermatrices}
t^{(n)} = tr_{V_0} R_{0 n} R_{0 n-1} \dots R_{01}
\end{equation}
with respect to the additional space $V_0$.
If the solution of the quantum Yang-Baxter equation depends on an additional spectral parameter
$\lambda$ then the corresponding transfer matrices commute:
\begin{equation}
\label{commut}
t^{(n)}(\lambda) t^{(n)}(\mu) = t^{(n)}(\mu) t^{(n)}(\lambda),
\end{equation}
which is the crucial fact in that theory (see \cite{B,TF,JM}).

In the general set-theoretical case we have a problem with the
trace operation (see although the remark at the end of this section), so one should look for an alternative. 
A possible candidate is given by the operators first considered by C.N Yang
in the classical paper \cite{Yang}. Their relation with
transfer matrices are nicely explained in Gaudin's book
\cite{Gaudin} (see especially Chapter 10, sections 2 and 3). A
related constructions were used by Frenkel and Reshetikhin in the
theory of q-KZ equation \cite{FR} and by Fomin and Kirillov in the
theory of Schubert polynomials \cite{FK} (see also \cite{BFZ, NY}).

Following \cite{V1} we introduce the {\it transfer maps} $T_i^{(n)}, i =
1,\dots, n$ as the maps of $X^{n}$ into itself defined by the
following formulas:
\begin{equation}
\label{T}
T_i^{(n)} = R_{i i+n-1} R_{i i+n-2} \dots R_{i i+1},
\end{equation}
where the indices are considered modulo $n$ with the agreement
that we are using $n$ rather than 0. In particular $T_1^{(n)} =
R_{1 n} R_{1 n-1} \dots R_{12}.$

{\bf Theorem \cite{V1}.} {\it For any reversible Yang-Baxter map}
$R$ {\it the transfer maps} $T_i^{(n)},\quad i =1,\dots, n$ {\it
commute with each other:}
\begin{equation}
\label{comm}
T_i^{(n)} T_j^{(n)} = T_j^{(n)} T_i^{(n)}
\end{equation}
{\it and satisfy the property}
\begin{equation}
\label{prod}
T_1^{(n)} T_2^{(n)} \dots T_n^{(n)} = Id.
\end{equation}
{\it Conversly, suppose that the maps} $T_i^{(n)}$ {\it determined by the formula (\ref{T}) commute
and satisfy the relation (\ref{prod}) for any $n \geq 2$ then $R$ is a reversible Yang-Baxter map.}

Proof of the first part is very similar to the proof of the commutativity of the transfer matrices
and follows from the consideration of the corresponding diagrams
representing the products $T_i^{(n)} T_j^{(n)}$ and $T_j^{(n)} T_i^{(n)}$ respectively (cf \cite{JM}):

\begin{figure}[htbp]
\begin{center}
\includegraphics[width=10cm]{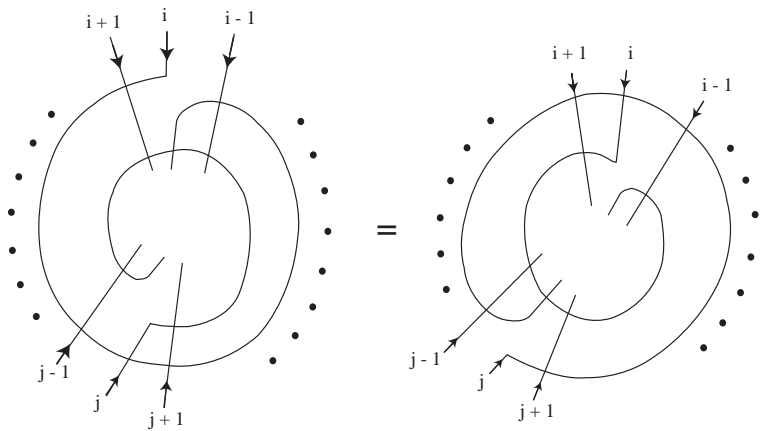}
\caption{Commutativity of the transfer maps}
\end{center}
\end{figure}



One can easily check that the second diagram is a result of the several operations presented
above on figures 1 and 2 applied to the first diagram.
The second identity also follows from similar consideration (see fig. below where the case $n=3$ is presented).

\begin{figure}[htbp]
\begin{center}
\includegraphics[width=10cm]{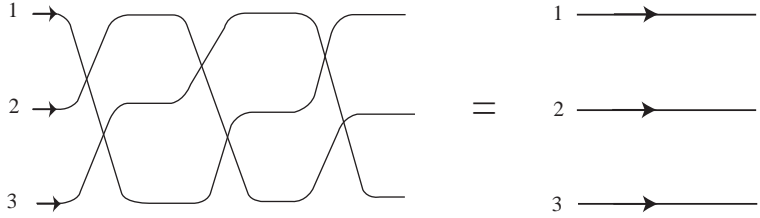}
\caption{Product of transfer maps}
\end{center}
\end{figure}


To prove the second part it is actually enough to consider only
the cases $n=2$ and $n=3.$ When $n=2$ we have two maps $T_1^{(2)}
= R$ and $T_2^{(2)} = R_{21}$, so the relation (\ref{prod})
becomes simply the reversibility condition (\ref{U}). The
commutativity condition is trivial in this case. For $n=3$ then we
have three transfer maps: $T_1^{(3)} = R_{13} R_{12}, T_2^{(3)} =
R_{21} R_{23}, T_3^{(3)} = R_{32} R_{31}.$ The product $T_1^{(3)}
T_2^{(3)} T_3^{(3)} =  R_{13} R_{12} R_{21} R_{23} R_{32} R_{31}$
is obviously an identity because $R_{12} R_{21}=R_{23}
R_{32}=R_{13} R_{31} =Id.$ Commutativity $T_1^{(3)} T_2^{(3)} =
T_2^{(3)} T_1^{(3)}$ means that $R_{13} R_{12} R_{21} R_{23} =
R_{21} R_{23} R_{13} R_{12}$ which is equivalent to Yang-Baxter
equation (\ref{YB}) modulo the reversibility relation which was
already shown.

The properties of the transfer maps can also be
explained in terms of the corresponding extended affine Weyl group
$\tilde A_{n-1}^{(1)}.$ As we have already mentioned above the
maps $$S_{i} = P_{ii+1} R_{ii+1},$$ where $P_{ij}$ is the
permutation of $i$-th and $j$-th factors in $X^{n},$ satisfy the
relations
\begin{equation}
\label{aff} S_{i}S_{i+1}S_{i} = S_{i+1}S_{i}S_{i+1}, \quad  i=1,
\dots, n-1
\end{equation}
and
\begin{equation}
\label{invol} S_i^2 = Id
\end{equation}
which are the defining relations of the affine Weyl group
$A_{n-1}^{(1)}.$ However (see the warning after the definition of
the Yang-Baxter maps) the transfer maps {\bf do not} belong to
this group but do belong to its extension $\tilde A_{n-1}^{(1)}$
generated by $S_1, \dots, S_n$ and the cyclic permutation $\omega
= P_{1 n} P_{1 n-1} \dots P_{12},\quad \omega^{n} = Id$ (known as
{\it extended affine Weyl group}, see \cite{Bour}). Moreover one
can easily check that they corresponds to the translations in this
groups and thus must commute (cf.\cite{NY}).

Note that if we consider only the set of $S_i$ with $i= 1, \dots,
n-1$ (without $S_n$) then they generate the (finite) permutation
group ${\bf S}_{n}.$ In that case as it was shown by Etingof et al
\cite{ESS} the corresponding dynamics is equivalent to the usual
action of ${\bf S}_{n}$ by permutations (and thus is trivial). It
is very important that in our case the transfer maps generate the
action of an infinite abelian group ${\bf Z}^{n-1}$ so the
question of whether two such actions are equivalent or not is a
non-trivial problem of dynamical systems theory. The fact that the
corresponding group is abelian indicates that the dynamics is
integrable but to make this claim precise is not so easy.

{\bf Remark. } We should mention that one can introduce another dynamical 
analogue of the transfer matrices,
which may be even more natural but leads to the multi-valued maps (correspondences). 
Indeed consider the map $$F=  R_{0 n} R_{0 n-1} \dots R_{01}: (x_0, x_1, \dots x_n) \rightarrow
(y_0, y_1, \dots, y_n)$$ and assume that $x_0 = y_0$ (to mimic the trace operation).
Generically one can express $x_0$ as a function of $ x_1, x_2, \dots, x_n $ from this relation
(in not necessarily unique way though). Substituting this $x_0$ into the rest of the relations 
$y_1 = y_1 (x_0, x_1, \dots, x_n), \dots,  y_n = y_n (x_0, x_1, \dots, x_n)$ we come to the {\it transfer correspondence}  $t^{(n)}: X^n \rightarrow X^n$, which can be considered as analogue of transfer matrices (\ref{transfermatrices}). The commutativity property (\ref{commut}) follows as in the usual case.
More detail can be found in \cite{GPV} where the corresponding multi-valued dynamics related to Adler maps and discrete KdV equation is investigated.

Now we are going to present the main constructions and the
examples of the Yang-Baxter maps coming from the theory of
integrable systems.


\section*{Matrix factorisations and Lax representations}

The relation between matrix factorisations and Yang-Baxter
property goes back at least to A.B. Zamolodchikov's works in the
late 70-th on factorised S-matrices and related algebras. I would
like to mention here also the papers \cite{ Skl, MV, R, HKKR, RS, O}. 
The first concrete example of the Yang-Baxter map I know
which came from this construction is the Adler map \cite{A}(see
below).

Let $A(x,\zeta)$ be a family of matrices depending on the point $x
\in X$ and a "spectral" parameter $\zeta \in {\bf C}$. One should
think of $X$ being an algebraic variety and $A$ depending
polynomially/rationally on $\zeta$ although the procedure works
always when the corresponding factorisation problems are uniquely
solvable. For matrix polynomials usually this is the case once the
factorisation of the corresponding determinant is fixed (see
\cite{GLR}).

The following procedure is a version of standard QR-algorithm in
linear algebra and known to be very useful in the theory of the
discrete integrable systems, see e.g. \cite{MV}.

Consider the product $L = A(y,\zeta) A(x,\zeta),$ then change the
order of the factors $L \rightarrow \tilde{L} = A(x, \zeta) A(y,
\zeta)$ and refactorise it again: $\tilde{L} = A(\tilde{y}, \zeta)
A(\tilde{x}, \zeta).$ Now define the map $R$ by the formula
\beq \label{fact} R(x,y) = (\tilde{x},\tilde{y}). \eeq 
In other words the relation
\begin{equation}\label{Lax}
 A(x,\lambda; \zeta)A(y,\mu; \zeta) =
 A(\tilde{y},\mu; \zeta)A(\tilde{x},\lambda; \zeta),
\end{equation}
is satisfied iff $(\tilde{x},\tilde{y}) = R (\lambda,\mu) (x,y).$

In that case the matrix $A(x,\zeta)$ is called the {\it Lax
matrix} of the map $R$, which is automatically reversible
Yang-Baxter map. Indeed the Yang-Baxter equation easily follows
from associativity and the reversibility is obvious.

Once the Lax matrix is known one can introduce the {\it monodromy
matrix} $$M= A(x_n,\zeta) A(x_{n-1},\zeta) \ldots A(x_1, \zeta).$$

{\bf Theorem \cite{V1}.} {\it The transfer maps} $T_i^{(n)}, \quad
i=1,\ldots, n$ {\it preserve the spectrum of the corresponding monodromy matrix}
$M(x_1, \ldots, x_n; \zeta)$ {\it for all} $\zeta.$

The proof is simple. For the map $T_1^{(n)} = R_{1 n} R_{1 n-1}
\dots R_{12}$ one can easily see that $A(x_n,\zeta)
A(x_{n-1},\zeta) \ldots A(x_1, \zeta) = A(\tilde{x}_1, \zeta)
A(\tilde{x}_n,\zeta) A(\tilde{x}_{n-1},\zeta)\ldots A(\tilde{x}_2,
\zeta),$ where $(\tilde{x}_1, \ldots ,\tilde{x}_n) = T_1^{(n)}
(x_1, \ldots, x_n).$ Since the spectrum of the product of the
matrices is invariant under the cyclic permutation of the factors,
we have $$Spec M(x_1, \ldots, x_n; \zeta) = Spec M(\tilde{x}_1,
\ldots ,\tilde{x}_n; \zeta).$$ To prove the same for $T_i^{(n)}$
one should consider the matrix $A(x_{n+i-1},\zeta) \ldots A(x_i,
\zeta)$ which obviously has the same spectrum as $M$ and use the
same arguments.

As a corollary we have that the coefficients of the characteristic
polynomial $\chi = \det (M - \lambda I)$ are the integrals of the
transfer maps $T_i^{(n)}.$ If the dependence of $A$ (and therefore
$M$) on $\lambda$ is polynomial then the dynamics can be
linearised on the Jacobi varieties of the corresponding spectral
curves $\chi(\zeta, \lambda) = 0$ (see \cite{MV}). If these maps
are also symplectic then one can use also the discrete version of
the Liouville theorem \cite{V} to claim their integrability.

We should mention here that although the procedure is similar to
the one proposed in \cite{MV} for the discrete Lagrangian systems
there is one important difference: in the construction of
Yang-Baxter maps the factors $A$ should depend only on one of the
variables, while in general scheme they usually depend both on $x$
and $y$.

{\bf Example.} Symmetries of the periodic dressing chain and
recutting of polygons \cite{VS, A}.

In \cite{A} Adler discovered a remarkable geometric dynamical
system, which he called the {\it recutting of polygons}. It turned
out that (after a slight modification) it can be described by the
following refactorisation procedure related to the {\it periodic
dressing chain} investigated by Shabat and the author \cite{VS}.

Here $X = {\bf C}\times{\bf C}$ and the matrix $A(x), x =
(f,\beta)$  has the form: 
\begin{equation}
\label{AdlerLax}
  A = \Bigl(
\begin{array}{cc}
f & f^2 + \beta - \zeta\\
1 & f
\end{array}
\Bigl) 
\end{equation} 

One can check that the factorisation procedure described
above leads in this case to the following birational map $R:
(f_1,\beta_1; f_2, \beta_2) \rightarrow (\tilde{f_1},\beta_1;
\tilde{f_2}, \beta_2):$
\begin{equation}
\label{Adler} \tilde{f_1} = f_2 - \frac{\beta_1 - \beta_2}{f_1 +
f_2} \qquad \tilde{f_2} = f_1 - \frac{\beta_2 - \beta_1}{f_1 +
f_2},
\end{equation}
which is the {\it Adler map}. Strictly speaking Adler considered
the $S$-version of this map: $S = P R$. The observation that it
satisfies the braid relations and thus leads to the action of the
affine Weyl group was important for him, because this was related
to the polynomial growth property \cite{V}. In  a different
context this action was discussed in much more details by Noumi
and Yamada \cite{NY}.

The transfer dynamics in this example is symplectic (see also
below) and integrable in the sense of discrete Liouville theorem
\cite{V} (see \cite{A}). This follows from the theory of the
dressing chain \cite{VS} which is closely related to the
finite-gap theory of the KdV equation \cite{Nov, DMN}. In fact the
transfer dynamics $T_i^{(n)}$ with odd $n=2g+1$ can be considered
as a discrete version of the Novikov equations \cite{Nov}
describing the $g$-gap Schr\"odinger operators (see for details
\cite{A,VS}).

It turned out that quite often (in particular, in the previous
example) one can derive the Lax matrix directly from the Yang-Baxter map. The
following observation is due to Suris and the author \cite{SV} and
was inspired by the results \cite{BS, N} on 3D consistent discrete
equations.

Suppose that on the set $X$ we have an effective action of the
linear group $G =  GL_N$, and that the Yang-Baxter map $R(\lambda,
\mu)$ has the following special form:
\begin{equation}\label{map}
\tilde{x} = B(y,\mu, \lambda)[x], \quad  \tilde{y} = A(x,\lambda,
\mu)[y],
\end{equation}
where $A, B: X \times {\bf C} \times {\bf C} \rightarrow GL_N$ are
some matrix valued functions on $X$ depending on parameters
$\lambda$ and $\mu$ and $A[x]$ denotes the action of $A \in G$ on
$x \in X.$ Then the claim is that both $A(x,\lambda, \zeta)$ and
$B^{\rm T}(x,\lambda,\zeta)$ are Lax matrices for $R$. The claim
about $B$ is equivalent to the fact that $B(x,\lambda,\zeta)$ is a
Lax matrix for $R_{21}.$

Indeed, look at the values of $z_{12}$ produced by the both parts
of the Yang-Baxter relation (\ref{sYB}): the left-hand side gives
$z_{12}=A(y_1,\mu,\nu)A(x_2,\lambda,\nu)[z]$, while the right-hand
side gives $z_{12}=A(x,\lambda,\nu) A(y,\mu,\nu)[z]$. Now since we
assume that the action of $G$ is effective, we immediately arrive
at the relation
$
A(x,\lambda,\nu) A(y,\mu,\nu)=A(y_1,\mu,\nu)A(x_2,\lambda,\nu),
$
which holds whenever $(x_2,y_1)=R(\lambda,\mu)(x,y)$. This
coincides with (\ref{fact}), an arbitrary parameter $\nu$ playing
the role of the spectral parameter $\zeta$.

Similarly, one could look at the values of $x_{23}$ produced by
the both parts of (\ref{sYB}): the left-hand side gives
$x_{23}=B(z,\nu,\lambda)B(y,\mu,\lambda)[x]$, while the right-hand
side gives $x_{23}=B(y_3,\mu,\lambda) B(z_2,\nu,\lambda)[x]$.
Effectiveness of the action of $G$ again implies:
$
B(z,\nu,\lambda)B(y,\mu,\lambda)=B(y_3,\mu,\lambda)
B(z_2,\nu,\lambda),
$
whenever $(y_3,z_2)=R(\mu,\nu)(y,z)$. This turns into (\ref{fact})
for the transposed matrices $B^{\rm T}$ (or for the inverse
matrices $B^{-1}$); the role of spectral parameter is here played
by an arbitrary parameter $\lambda$.

In fact the effectiveness of the action can be replaced by the
following weaker property. Let us call the action of $G = GL_N$ on
$X$ {\it projective} if the action of $A$ on $X$ is trivial only
if $A$ is a scalar. Similarly we say that $A(x, \lambda,\zeta)$
gives a {\it projective Lax representation} for the Yang-Baxter
map $R$ if the relation (\ref{Lax})holds up to multiplication by a
scalar matrix $c I$, where $c$ may depend on all the variables in
the relation. One can easily modify the previous arguments to
produce the integrals for the transfer maps using the projective
Lax matrix: all the ratios of the eigenvalues of the monodromy
matrix are obviously preserved by these maps.

For the projective action the matrices $A(x,\lambda,\zeta)$ and
$B^{\rm T}(x,\lambda,\zeta)$ give projective Lax representations
for the corresponding Yang-Baxter maps (\ref{fact}). In practice
however for a natural choice of matrices $A$, $B$ we have actually
proper Lax representations, as we can see from the example of
Adler map (\ref{Adler}).

Indeed one can rewrite (\ref{Adler}) as follows $$\tilde{y} = x -
\frac{\mu - \lambda}{x+y} = \frac{x^2 + xy -(\mu - \lambda)}{x+y}
= A(x, \lambda, \mu)[y],$$ where as before
 $$ A(x,\lambda,\zeta) = \Bigl(\begin{array}{cc}
 x &  x^2 + \lambda - \zeta\\ 1 & x \end{array} \Bigl). $$
and the group $G = GL_2$ is acting on ${\bf CP}^1$ by M\"obius
transformations.

Another interesting example when the Lax representation can be
read directly from the form of the Yang-Baxter map is related to
the theory of solitons, in particular matrix solitons.

\section*{Interaction of solitons as a Yang-Baxter map.}

Consider any integrable by the inverse scattering method PDE in
1+1 dimensions which has multisoliton solutions. Suppose that each
soliton has a non-trivial internal degrees of freedom described by
the set $X$ (which is usually a manifold). The soliton interaction has the following 
remarkable property: the  result 
of interaction of three solitons is independent of the order of collisions
and thus is completely defined by the pairwise interaction.
Although the proof of this claim in full generality seems to be absent in the literature
it is confirmed in all known cases (see recent papers  \cite{G, Tsu, APT} 
for the discussion of this). 
This property means that the corresponding soliton interaction map $R$  satisfies
the Yang-Baxter relation.

A good example is the matrix KdV equation \cite{Lax} : $$
U_t+3UU_x+3U_xU+U_{xxx}=0, $$ $U$ is $d \times d$ matrix.

It is easy to check that it has the soliton solution of the form
$$U = 2 \lambda^2 P \sech^2 (\lambda x- 4\lambda^3 t),$$ where $P$
must be a projector: $P^2 = P.$ If $P$ has rank 1 then it should
has the form $P=\frac{\dpl \xi\otimes \eta}{\dpl (\xi,\eta)},$
where $\xi$ is a vector in a vector space $V$ of dimension $d$,
$\eta$ is a vector from the dual space $V^*$ (covector) and
bracket $(\xi,\eta)$ means the canonical pairing between $V$ and
$V^*.$

The change of the matrix amplitudes $P$ ("polarisations") of two
solitons with the velocities $\lambda_1$ and $\lambda_2$ after
their interaction is described by the following map \cite{G,GV}:
$$R(\lambda_1, \lambda_2): (\xi_1, \eta_1; \xi_2, \eta_2)
\rightarrow (\tilde{\xi_1}, \tilde{\eta_1}; \tilde{\xi_2},
\tilde{\eta_2)}$$ \beq \label{fi} \tilde{\xi_1} = \xi_1+\frac{\dpl
2\lambda_2(\xi_1,\eta_2)}{\dpl (\lambda_1-
\lambda_2)(\xi_2,\eta_2)}\xi_2, \qquad \tilde{\eta_1} =
\eta_{1}+\frac{\dpl 2\lambda_2(\xi_2,\eta_1)}{\dpl (\lambda_1-
\lambda_2)(\xi_2,\eta_2)}\eta_2, \eeq

\beq \label{fj} \tilde{\xi_2} = \xi_2+\frac{\dpl
2\lambda_1(\xi_2,\eta_1)}{\dpl (\lambda_2-
\lambda_1)(\xi_1,\eta_1)}\xi_1,\qquad \tilde{\eta_2} =
\eta_2+\frac{\dpl 2\lambda_1(\xi_1,\eta_2)}{\dpl (\lambda_2
-\lambda_1)(\xi_1,\eta_1)}\eta_1. \eeq

In this example $X$ is the set of projectors $P$ of rank 1 which
is the variety ${\bf CP}^{N-1} \times {\bf CP}^{N-1}$, and the
group $G = GL_N$ is acting on the projectors by conjugation (which
corresponds to the natural action of $GL(V)$ on $V \otimes V^*$).

It is easy to see that the formulas (\ref{fi}), (\ref{fj}) are of
the form (\ref{map}) with the matrices
\begin{equation}
\label{Lax1} A(\xi, \eta, \lambda; \zeta) = I + \frac{2
\lambda}{\zeta -\lambda}\frac{\dpl \xi\otimes \eta}{\dpl
(\xi,\eta)}= I+\frac{2\lambda}{\zeta-\lambda}P
\end{equation}
(note that again $R_{21}=R$). Our results show that the matrix
$A(P,\lambda,\zeta)$ gives a projective Lax representation for the
interaction map.

One can show that this is actually a genuine Lax representation.
We will demonstrate this for the class of more general Yang-Baxter
maps on Grassmannians following to \cite{GV}.

Let $V$ be an $n$-dimensional real (or complex) vector space, $P:
V \rightarrow V$ be a projector  of rank $k$: $P^2 =P$. Any such
projector is uniquely determined by its kernel $K = Ker P$ and
image $L = Im P,$ which are two subspaces of $V$ of dimensions $k$
and $n-k$ complementary to each other: $K \oplus L = V.$ The space
of all projectors $X$ of rank $k$ is an open set in the product of
two Grassmannians $G(k,n) \times G(n-k,n).$

Consider the following matrix
\begin{equation}
\label{Lax2} A(P, \lambda; \zeta) = I + \frac{2 \lambda}{\zeta
-\lambda} P
\end{equation}
and the related re-factorization relation
\begin{equation}
\label{rat}
 (I + \frac{2 \lambda_1}{\zeta -\lambda_1} P_1)(I +
\frac{2 \lambda_2}{\zeta -\lambda_2} P_2) = (I + \frac{2
\lambda_2}{\zeta -\lambda_2} \tilde{P}_2)(I + \frac{2
\lambda_1}{\zeta -\lambda_1} \tilde{P}_1)
\end{equation}
which we can rewrite in the polynomial form as
\begin{equation}
\label{ref} ((\zeta -\lambda_1) I + 2 \lambda_1 P_1)((\zeta
-\lambda_2) I + 2 \lambda_2 P_2) = ((\zeta -\lambda_2) I + 2
\lambda_2\tilde{P}_2)((\zeta -\lambda_1) I + 2 \lambda_1
\tilde{P}_1).
\end{equation}

The claim is that if $\lambda_1 \neq \pm\lambda_2$ it has a unique
solution. This follows from the general theory of matrix
polynomials (see e.g. \cite{GLR}) but in this case we can see this
directly.

Indeed let us compare the kernels of both sides of the relation
(\ref{ref}) when the spectral parameter $\zeta = \lambda_1.$ In
the right hand side we obviously have $\tilde{K}_1$ while the left
hand side gives $$((\lambda_1 -\lambda_2) I + 2 \lambda_2
P_2)^{-1} K_1 = (I + \frac{2 \lambda_2}{\lambda_1
-\lambda_2}P_2)^{-1} K_1.$$

Now we use the following property of the matrix (\ref{Lax2}):
\begin{equation}
\label{property} A(P, -\lambda; \zeta) =  A(P, \lambda;
\zeta)^{-1}
\end{equation}
to have
\begin{equation}
\label{K1}\tilde{K}_1 = (I - \frac{2 \lambda_2}{\lambda_1
+\lambda_2}P_2) K_1.
\end{equation}

Similarly taking the image of both sides of (\ref{ref}) at $\zeta
= \lambda_2$ we will have
\begin{equation}
\label{L2}\tilde{L}_2 = (I + \frac{2 \lambda_1}{\lambda_2
-\lambda_1}P_1) L_2.
\end{equation}

To find $\tilde{K}_2$ and $\tilde{L}_1$ one should take first the
inverse of both sides of (\ref{rat}), use the property
(\ref{property}) and then repeat the procedure. This will lead us
to the formulas:
\begin{equation}
\label{K2}\tilde{K}_2 = (I - \frac{2 \lambda_1}{\lambda_1
+\lambda_2}P_1) K_2
\end{equation}
and
\begin{equation}
\label{L1}\tilde{L}_1 = (I + \frac{2 \lambda_2}{\lambda_1
-\lambda_2}P_2) L_1.
\end{equation}

The formulas (\ref{K1},\ref{L2},\ref{K2},\ref{L1}) determine a
parameter-dependent Yang-Baxter map on the set of projectors. One
can easily check that for $k=1$ one has the formulas (\ref{fi},
\ref{fj}) for two matrix soliton interaction.

If we supply now our vector space $V$ with the Euclidean
(Hermitian) structure and consider the self-adjoint projectors $P$
of rank $k$ then the corresponding space $X$ will coincide with
the Grassmannian $G(k,n):$ such a projector is completely
determined by its image $L$ (which is a $k$-dimensional subspace
in $V$ and thus a point in $G(k,n)$) since the kernel $K$ in this
case is the orthogonal complement to $L.$

The corresponding Yang-Baxter map $R$ on the Grassmannian is
determined by the formulas
\begin{equation}
\label{L11}\tilde{L}_1 = (I + \frac{2 \lambda_2}{\lambda_1
-\lambda_2}P_2) L_1,
\end{equation}
\begin{equation}
\label{L22}\tilde{L}_2 = (I + \frac{2 \lambda_1}{\lambda_2
-\lambda_1}P_1) L_2.
\end{equation}


\subsection*{Yang-Baxter maps related to geometric crystals} 

In this section we briefly discuss the examples of Yang-Baxter maps related to geometric crystals 
and tropical combinatorics \cite{BK, E, KNY, KOTY}. We show how one can derive the Lax representation for them following to \cite{SV}.

Let $X={\bf C}^n$, and define the map $R:X\times X\rightarrow X\times X$
by the formulas
\begin{equation}\label{E}
 \tilde{x}_j=x_j\,\frac{P_j}{P_{j-1}},\qquad
 \tilde{y}_j=y_j\,\frac{P_{j-1}}{P_j},\qquad j=1,\ldots,n,
\end{equation}
 where
\begin{equation}\label{EP}
 P_j=\sum_{a=1}^n\left(\prod_{k=1}^{a-1}x_{j+k}\prod_{k=a+1}^ny_{j+k}\right)
\end{equation}
(in this formula subscripts $j+k$ are taken (mod $n$)).
The map (\ref{E}) obviously preserves the following subsets invariant:
 $X_\lambda\times X_\mu\subset X\times X$, where
 $X_\lambda=\{(x_1,\ldots,x_n)\in X: \prod_{k=1}^nx_k=\lambda\}$.
 Let us show that the restriction of $R$ to $X_\lambda\times X_\mu$
 can be written in the form (\ref{map}). For this,  let us embed this set into
 ${\bf CP}^{n-1}\times{\bf CP}^{n-1}$:
\[
J(x,y)=(z(x),w(y)),\quad z(x)=(1:z_1:\ldots:z_{n-1}), \quad
w(y)=(w_1:\ldots:w_{n-1}:1),
\]
\[
z_j=\prod_{k=1}^{j}x_k\,,\qquad w_j=\prod_{k=j+1}^ny_k\,.
\]
Then it is easy to see that in coordinates $(z,w)$ the map $R$ is
written as
\[
\tilde{z}=B(y,\mu,\lambda)[z]\,,\qquad
\tilde{w}=A(x,\lambda,\mu)[w]\,,
\]
with certain matrices $B,A$ from $G=GL_n$, where the standard
projective action of $GL_n$ on ${\bf CP}^{n-1}$ is used. 
A simple calculation shows that the inverse matrices are the following
matrices well-known in this context \cite{E,KNY, KOTY}:
\begin{equation}\label{B}
B^{-1}(y,\mu,\lambda)=\left(\begin{array}{cccccc}
 y_1 & -1 & 0 & \ldots & 0 & 0\\
 0 & y_2 & -1 & \ldots & 0 & 0\\
 0 & 0 & y_3 & \ldots & 0 & 0 \\
 & \ldots & &&\ldots &\\
 0 & 0 & 0 & \ldots & y_{n-1} & -1\\
 -\lambda &  0 & 0 & \ldots & 0 & y_n
 \end{array}\right),
\end{equation}
\begin{equation}\label{A}
A^{-1}(x,\lambda,\mu)=\left(\begin{array}{cccccc}
 x_1 & 0 & 0 & \ldots & 0 & -\mu\\
 -1 & x_2 & 0 & \ldots & 0 & 0\\
 0 & -1 & x_3 & \ldots & 0 & 0 \\
 & \ldots & &&\ldots &\\
 0 & 0 & 0 & \ldots & x_{n-1} & 0\\
 0 &  0 & 0 & \ldots & -1 & x_n
 \end{array}\right).
\end{equation}
Thus we derived again the Lax matrices directly from the map itself.
It is worthy to mention that this is closely related to the notion of the {\it
structure group} $G_R$ of the Yang-Baxter map $R$ \cite{E}. It
was shown by Etingof in \cite{E} that for the map (\ref{E}) the
so-called reduced structure groups $G_R^+$ and $G_R^-$ can be
realised as the subgroups of the loop group $PGL_n(C(\lambda))$
generated by the matrix functions $A^{-1}(x,\cdot,\lambda)$ with $x\in X$, 
resp. by $B^{-1}(x,\cdot,\lambda)$ with $x\in X$.

\section*{Poisson Lie groups and Hamiltonian structure for the transfer dynamics}

In this section we describe the  Hamiltonian structure for the transfer maps following 
the recent work by Reshetikhin and the author \cite{RV} based on some general constructions from
Weinstein and Xu \cite{WX}. The main observation is that if the Lax matrices form symplectic leaves in certain Poisson Lie group then the transfer dynamics is Poisson.

For a good account of the theory of the Poisson Lie group we
recommend Chapter 5 in the Reiman and Semenov-Tian-Shansky book
\cite{RSTS} or original Drinfeld paper \cite{Drin}.

Recall that {\it Poisson Lie group} is a Lie group $G$ with the
Poisson structure on it such that the multiplication map $\mu: G
\times G \rightarrow G$ is Poisson. Recall also that the map
between two Poisson manifolds is called {\it Poisson} if it
preserves the Poisson bracket.

We will need also the notion of {\it Poisson correspondence} (or
Poisson relation) (see e.g. Weinstein \cite{W} and references therein). The
{\it correspondence} (relation) $\Phi$ from a set $X$ into itself
is a multivalued, partially defined map determined by its graph
$\Gamma_{\Phi}$, which is a subset of $X \times X.$ Now let $X =
M$ be a smooth manifold with Poisson structure $J$ and the graph
$\Gamma_{\Phi}$ is a submanifold of $M \times M.$ Let us supply
the manifold $M \times M$ with the Poisson structure $J \oplus (-
J).$ Recall now that a submanifold $N$ of a Poisson manifold $P$
is called { \it coisotropic} if the set of functions on $P$ which
vanish on $N$ is closed under Poisson bracket. The correspondence
$\Phi$ is called {\it Poisson} if its graph $\Gamma_{\Phi}$ is
coisotropic submanifold of $M \times M.$ For the usual maps this
equivalent to the standard notion of Poisson map. If some
smoothness condition is fulfilled then the composition of two
Poisson correspondences is again Poisson.

For any Poisson Lie group $G$ one can construct the following
correspondence $\Phi_G: G \times G \rightarrow G \times G$: its
graph $\Gamma_{\Phi}$ is the set of  $(g_1, g_2; h_1, h_2) \in G
\times G \times G \times G$ such that
\begin{equation}
g_1 g_2 = h_2 h_1. \label{factor}
\end{equation}

From the definition of the Poisson Lie group it follows that this
correspondence is Poisson.

The first simple observation is that this correspondence satisfies
the Yang-Baxter relation. Note that for the correspondence $R: X
\times X \rightarrow X \times X$ the map $R_{12}$ is defined by
its graph which is a subset of $X^{(6)}$ of the form $\Gamma_R
\times \Delta,$ where $\Delta = (x,x), x \in X$ is the diagonal in
$X \times X.$ Similarly one can define $R_{ij}$ and thus {\it
Yang-Baxter correspondence} as such $R$ which satisfies the
relation (\ref{YB}). Reversibility is a bit more tricky. The
problem is with the notion of the inverse correspondence
$\Phi^{-1}$: the standard formula $\Phi \Phi^{-1} = Id$ seems to
be too restrictive, so we should avoid it. There are two
involutions on $X \times X \times X \times X$: $$ \pi:  (x_1, x_2;
y_1, y_2) = (x_2, x_1; y_2, y_1)$$ and $$\tau:  (x_1, x_2; y_1,
y_2) =  (y_1, y_2; x_1, x_2).$$ We will call the correspondence
$R$ {\it reversible} if $\pi \Gamma_{R} = \tau \Gamma_{R},$ or
equivalently its graph $\Gamma_{R}$ is invariant under the
involution $$\sigma = \pi \tau: (x_1, x_2; y_1, y_2) = (y_2, y_1;
x_2, x_1).$$ For the maps this is equivalent to the definition
given above.

We claim that for any Poisson Lie group $G$ the
Poisson correspondence $\Phi_G$ is a reversible Yang-Baxter
correspondence. Indeed the Yang-Baxter property  follows from the associativity 
of multiplication in $G$: both left and right sides are the correspondences with the same graph given by the relation $$g_1 g_2 g_3 = h_3 h_2 h_1.$$ Reversibility is obvious.

In such generality this claim is probably not interesting enough,
since we are primarily interested in genuine maps but not
correspondences.  A natural idea therefore is to restrict the
previous construction to the symplectic leaves in $G$ in order to
make the re-factorisation relation (\ref{factor}) uniquely solvable.

Suppose that there exist a one-parameter family of embeddings of
the set $X$ as a symplectic leaf in the Poisson Lie group $G$:
$\phi_{\lambda} : X \rightarrow G$ and define the correspondence
$R(\lambda, \mu): X \times X \rightarrow X \times X$ by the
relation
\begin{equation}
\phi_{\lambda}(x) \phi_{\mu}(y) =  \phi_{\mu}(\tilde{y})
\phi_{\lambda}(\tilde{x}). \label{factl}
\end{equation}

Note that the embedding $\phi_{\lambda} $ induce on $X$ the
symplectic structure
\begin{equation}
\label{symp} \omega_{\lambda} = \phi_{\lambda}^ *(\omega_{G}),
\end{equation}
where $\omega_{G}$ is the symplectic structure on the
corresponding leaf in $G$. Let us introduce on $X \times X$ a
symplectic structure as the direct sum $\omega_{\lambda} \oplus
\omega_{\mu}.$ In a similar way we define the symplectic structure
$\Omega^{(N)}$ on $X^{(N)}$ as the direct sum of the forms
$\omega_i = \omega_{\lambda_i}:$
\begin{equation}
\label{Omega} \Omega^{(N)} = \omega_{1} \oplus \omega_{2} \oplus
...\oplus \omega_{N}.
\end{equation}

{\bf Theorem \cite{WX, RV}.} {\it  The correspondence $R(\lambda,
\mu)$ defined by (\ref{factl})  is a reversible Yang-Baxter
Poisson correspondence. The related transfer maps (\ref{T}) are
Poisson correspondences with respect to the symplectic structure}
$\Omega^{(N)}.$

This construction is already good enough to produce many
interesting examples of symplectic Yang-Baxter maps, including the
interaction of matrix solitons \cite{V1,GV}. To show this we
consider the standard Poisson Lie loop group $G = LGL(n, \bf{R})$
and the symplectic leaves of a special form.

The standard Poisson structure on the loop group $G = LGL(n,
\bf{R})$ is given by the Sklyanin formula
\begin{equation}
\{L_1(\zeta),\otimes L_2(\eta) \} = [ \frac{P_{12}}{\zeta - \eta},
L_1(\zeta) \otimes L_2(\eta)], \label{Skl}
\end{equation}
where $P_{12}$ is the permutation matrix: $P(x \otimes y) = y
\otimes x.$

We will consider the symplectic leaves of the form $L(\zeta) = A +
B \zeta.$ The relation (\ref{Skl}) is equivalent to the following
set of relations
\begin{equation}
\label{leaves} \{A_1, \otimes A_2\} = P_{12}(A_1 B_2 - B_1 A_2),
\quad \{A_1, \otimes B_2\} = 0, \quad \{B_1, \otimes B_2\} = 0.
\end{equation}

From these relations it follows that the matrix elements of $B$
are in the centre of the corresponding Poisson algebra, so $B$ can
be fixed.

If $B$ is non-degenerate one can assume that $B = I$ is identity
matrix. The rest of the relations (\ref{leaves}) determine the
standard Lie-Poisson structure on the matrix group $GL(n, \bf{R})$
considered as an open set in the $n \times n$-matrix algebra $M_n$
which can be identified as a space with $gl^*(n, \bf{R})$.

Let ${\mathcal O}_k$ be the special symplectic leaf in $gl^*(n,
\bf{R})$ consisting of involutions $S = 2P - I,$ where $P$ is a
projector of rank $k$: $P^2 = P, S^2 = I.$ Such a projector $P$ is
determined by its kernel $K = Ker P$ and its image $L = Im P,$
which are two subspaces of dimension $k$ and $n-k$ respectively
such that $K \oplus L = \bf{R}^n.$ Thus $X$ is an open subset of
the product of two Grassmannians $G(n, k)$ and $G(n, n-k)$. Note
that the Grassmannian $G(n, n-k)$ is naturally isomorphic to the
Grassmannian $G(n, k)$ in the dual space.

Let $X = {\mathcal O}_k$ and define the embedding $\phi_{\lambda}$
of $X$ into $G$ by the formula
\begin{equation}
\label{phi} \phi_{\lambda}(S) = \zeta I + \lambda S = (\zeta -
\lambda) I + 2 \lambda P .
\end{equation}

It is easy to see that the formulas (\ref{factl}), (\ref{phi}) lead to the Yang-Baxter map 
(\ref{K1},\ref{L2},\ref{K2},\ref{L1}) considered in the previous section in  relation with matrix
KdV equation.  In the rank 1 case we have the formulas  (\ref{fi}), (\ref{fj}). 

The corresponding transfer maps $T_i^{(N)} $ (\ref{T}) are
symplectic and preserve the spectrum of the monodromy matrix
\begin{equation}
\label{mon} M(\zeta) =  {\mathcal L}(P_1,\lambda_1, \zeta)
{\mathcal L}(P_{2},\lambda_{2},\zeta) \ldots {\mathcal
L}(P_N,\lambda_N, \zeta),
\end{equation}
where ${\mathcal L}(P, \lambda, \zeta) = (\zeta - \lambda) I + 2
\lambda P.$ The coefficients of the characteristic polynomial
$\chi(\zeta, \eta) = \det ({\mathcal L} - \eta I)$ give the set of
integrals for transfer dynamics, which are in involution with
respect to the symplectic structure $\Omega^{(N)}$ on $X^{(N)}.$

We believe that this implies integrability of the transfer dynamics
in Liouville sense, but the proof of completeness is still to be found. 
As usual the dynamics is linear on the
Jacobi variety of the spectral curve $\chi(\zeta, \eta) = 0.$

Let us consider now the symplectic leaves in  $LGL(2, {\bf
R})$ of the form $L = A + B \zeta$ with degenerate $$ B = \Bigl(
\begin{array}{cc}
0 &  1 \\ 0  & 0
\end{array}
\Bigl). $$ Writing $A$ as $$ A = \Bigl(
\begin{array}{cc}
a &  b \\ c  & d
\end{array}
\Bigl) $$ we have the following Poisson brackets for the matrix
elements from the relation (\ref{Skl}): $$ \label{abcd} \{a, b\} =
a, \quad \{a, d\} = c, \quad \{a, c\} = 0, \quad \{b, c\} = 0,
\quad \{b, d\} = d, \quad \{c, d\} = 0. $$ It is easy to check
that the centre of this Poisson algebra is generated by two
functions: $C_1 = c$ and $C_2 = ad - bc.$ If we choose $C_1 = c
=1$ and $C_2 = ad - bc = \lambda$ we come to the
symplectic leaf ${\mathcal O}(\lambda)$ consisting of the matrices of
the form
\begin{equation}
\label{Leaf}
 L (\zeta, \lambda) =
\Bigl(
\begin{array}{cc}
a &  ad + \lambda - \zeta\\ 1 & d
\end{array}
\Bigl),
\end{equation}
where $\{a, d\} = 1.$

This is not quite Lax matrix (\ref{AdlerLax}) for the Adler map yet. To have it one should apply the Hamiltonian reduction to the product of the corresponding orbits ${\mathcal O}(\lambda)$ (see \cite{RV} for the details). This leads to the Poisson structure for the transfer dynamics for Adler maps, which is known from the theory of the periodic dressing chain \cite{VS, A} (with an odd period $N$).  

\section*{Classification of Yang-Baxter maps: quadrirational case}

The classification of the Yang-Baxter maps is probably the most interesting but hopelessly difficult problem even in the simplest case of $X={\bf CP}^1$. A natural class where the complete classification seems to be realistic is given by the quadrirational maps introduced and investigated by  
 Adler, Bobenko and Suris \cite{ABS2}. 
 
The map $R: (x,y) \rightarrow (u,v)$ is called {\it quadrirational} if both $R$ and so called {\it companion map} 
$\bar R: (x, v) \rightarrow (y,u)$ are birational maps of ${\bf CP}^1\times{\bf CP}^1$ into itself.
All such maps have the form
\begin{equation}\label{maps}
 R:\;\; u=\frac{a(y)x+b(y)}{c(y)x+d(y)}\,,\quad
         v=\frac{A(x)y+B(x)}{C(x)y+D(x)}\,,
\end{equation}
where $a(y),\ldots,d(y)$ and $A(x),\ldots,D(x)$ are
polynomials of degree at most 2. There are three subclasses of such
maps, denoted as [1:1], [1:2] and [2:2], corresponding to the highest
degrees of the coefficients of both fractions in (\ref{maps}). 
In \cite{ABS2} the classification of such maps was done modulo group $({\cal M}\ddot{o}b)^4$, where
${\cal M}\ddot{o}b$ is the group of projective transformations of ${\bf CP}^1$:
$$z \rightarrow w = \frac{az+b}{cz+d}.$$ In the richest class [2:2] it looks as follows.

{\bf Theorem \cite{ABS2}.}
{\it Any quadrirational map on ${\bf CP}^1\times{\bf CP}^1$ of subclass [2:2] is equivalent,
under some change of variables from $({\cal M}\ddot{o}b)^4$, to exactly one
of the following five maps:}
\begin{align}
\label{F1}\tag{$F_\I$}
  u&= \a yP,&
  v&= \b xP,&
  P&= \frac{(1-\b)x+\b-\a+(\a-1)y}
           {\b(1-\a)x+(\a-\b)yx+\a(\b-1)y},& \\
\label{F2}\tag{$F_\II$}
  u&= \frac{y}{\a}\,P,&
  v&= \frac{x}{\b}\,P,&
  P&= \frac{\a x-\b y+\b-\a}{x-y},& \\
\label{F3}\tag{$F_\III$}
  u&= \frac{y}{\a}\,P,&
  v&= \frac{x}{\b}\,P,&
  P&= \frac{\a x-\b y}{x-y},& \\
\label{F4}\tag{$F_\IV$}
  u&= yP,&
  v&= xP,&
  P&= 1+\frac{\b-\a}{x-y},& \\
\label{F5}\tag{$F_\V$}
  u&= y+P,&
  v&= x+P,&
  P&= \frac{\a-\b}{x-y},&
\end{align}
{\it with some suitable constants} $\a,\b$.

It is remarkable that all these maps can be described by the following natural geometric construction \cite{ABS2}.

\begin{figure}[htbp]
\begin{center}\includegraphics[width=8cm]{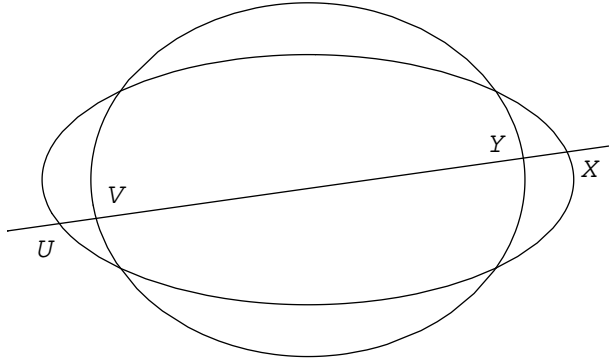}\end{center}
\caption{A quadrirational map on a pair of conics}\label{fig:greeks}
\end{figure}

Consider a pair of nondegenerate conics $Q_1$, $Q_2$ on the plane ${\bf CP}^2$, which are rational curves isomorphic to ${\bf CP}^1$.  Take $X\in Q_1$,
$Y\in Q_2$, and let $\ell=\overline{XY}$ be the line through $X,Y$. 
Generically the line $\ell$ intersects  $Q_1$ and $Q_2$ at
two points $U$ and $V$ (see Fig.~\ref{fig:greeks}). 
This defines a map $\cR:(X,Y)\mapsto(U,V)$. 

There are five possible types $\I-\V$ of intersection of two conics (see  \cite{Berger}):
\begin{itemize}\setlength{\itemsep}{0pt}
\item[\I:] four simple intersection points;
\item[\II:] two simple intersection points and one point of tangency;
\item[\III:] two points of tangency;
\item[\IV:] one simple intersection point and one point of the second order
tangency;
\item[\V:] one point of the third order tangency.
\end{itemize}
In suitable rational parametrisations of the conics
we have exactly the five maps listed in the theorem above.


It is remarkable that all five representative maps from the last theorem satisfy the Yang-Baxter relation (see \cite{ABS2}, where a nice geometric interpretation of this relation is also given). This surprising fact, which must be related to the classification approach of \cite{ABS2} based on the analysis of the singularities,  does not imply however that all quadrirational maps are Yang-Baxter (as one might conclude from \cite{ABS2}) since the action of $({\cal M}\ddot{o}b)^4$ in general destroy the Yang-Baxter property. 

For example, if we change in the fifth map $F_V$ the variables $x \rightarrow -x, y \rightarrow -y$ (preserving $u,v$)
we come to the quadrirational map  $u = -y- P,  v= -x- P, 
P= \frac{\a-\b}{x-y}, $ which does not satisfy the Yang-Baxter relation.
However, if we change $u \rightarrow -u, y \rightarrow -y$ we come to the Adler map (\ref{Adler}),
 which is a Yang-Baxter map. Note that neither of these changes is a proper conjugation of maps
 of ${\bf CP}^1\times{\bf CP}^1$ into itself (in particular, $F_V$ is not conjugated to Adler map).

A question remains therefore to classify all quadrirational Yang-Baxter maps modulo diagonal action of ${\cal M}\ddot{o}b$, which seems to be a natural equivalence within this class. The previous example shows that in the same $({\cal M}\ddot{o}b)^4$-orbit there are several Yang-Baxter representatives, non-equivalent under ${\cal M}\ddot{o}b$-conjugation. The answer to this question can be probably extracted from the results of the paper \cite{ABS2}.

Let us mention in this relation that a complete classification of integrable (3D-consistent) discrete equations on quad-graphs is known due to the same authors \cite{ABS}. Since these two problems are closely related (see \cite{ABS, PTV}) a very interesting question is what could be a Yang-Baxter analogue of the "master equation" $Q_4$ in their classification (see \cite{AS} for the discussion of this very important case). Probably one should look for a wider class of Yang-Baxter maps (like correspondences or Cremona maps). A class of matrix factorisations suggested by Odesskii  \cite{O} could give one of the possible answers to this question.

\section*{ Acknowledgements}

I am very grateful to the organisers of the meeting "Combinatorial
Aspects of Integrable Systems" (Kyoto, July 2004), for inviting me
to give these lectures.

I am also very grateful to V. Adler, V. Drinfeld, P. Etingof, S. Fomin, V.
Papageorgiou, N. Reshetikhin and Yu. Suris for useful discussions and comments.

This work was partially supported by the European network ENIGMA
(contract MRTN-CT-2004-5652) and the ESF programme "Methods of Integrable Systems, Geometry, Applied Mathematics" (MISGAM).


\begin{thebibliography}{99}
\bibitem{D}
V.G. Drinfeld {\it On some unsolved problems in quantum group theory.}
In "Quantum groups" (Leningrad, 1990), Lecture Notes in Math., 1510, Springer, 1992, p. 1-8.
\bibitem{Skl}
E.K. Sklyanin {\it Classical limits of ${\rm SU}(2)$-invariant solutions
of the Yang-Baxter equation.} J. Soviet Math. 40 (1988), no. 1,
93--107.
\bibitem{WX}
A. Weinstein, P. Xu {\it Classical solutions of the quantum Yang-Baxter equation.}
Comm. Math. Phys. 148 (1992), 309-343.
\bibitem{ESS}
P. Etingof, T. Schedler, A. Soloviev {\it Set-theoretical solutions to the quantum Yang-Baxter equation.}
Duke Math. J., V. 100 (1999), 169-209.
\bibitem{Solo}
A. Soloviev {\it Non-unitary set-theoretical solutions to the quantum Yang-Baxter equation.}  Math. Res. Lett. {\bf 7} (2000), no. 5-6, 577--596.
\bibitem{Hi}
J. Hietarinta {\it Permutation-type solutions to the Yang-Baxter and other n-simplex equations.}
J. Phys. A: Math. Gen. 30 (1997), 4757-4771.
\bibitem{Buch}
V.M. Buchstaber {\it The Yang-Baxter transformation.} Russ. Math. Surveys, 53:6 (1998), 1343-1345.
\bibitem{ESG}
P. Etingof, A. Soloviev, R. Guralnick {\it  Indecomposable set-theoretical solutions to the quantum Yang-Baxter equation on a set with a prime number of elements. } J. Algebra {\bf 242} (2001), no. 2, 709--719.
\bibitem{LYZ}
J.-H. Lu, M. Yan, Y.-C.  Zhu
{\it On the set-theoretical Yang-Baxter equation.}
Duke Math. J. {\bf 104} (2000), no. 1, 1--18.
\bibitem{GB}
P. Gu, C.-M. Bai {\it On set-theoretical solutions to quantum Yang-Baxter equation.}  Commun. Theor. Phys. (Beijing) {\bf 39} (2003), no. 2, 141--143.
\bibitem{G-I}
T. Gateva-Ivanova
{\it A combinatorial approach to the set-theoretic solutions of the  Yang-Baxter equation.} (2004) math.QA/0404461. J. Math. Phys. {\bf 45} (2004), no. 10, 3828--3858.
\bibitem{CS}
J. Scott Carter, M. Saito {\it Set-theoretical Yang-Baxter solutions via Fox calculus.} (2005)
GT/0503166

\bibitem{BK}
A. Berenstein, D. Kazhdan {\it Geometric and unipotent crystals.}
GAFA 2000 (Tel Aviv, 1999).  Geom. Funct. Anal. 2000, Special Volume, Part I, 188--236
\bibitem{E}
P. Etingof {\it Geometric crystals and set-theoretical solutions to the quantum Yang-Baxter equation.}
Comm. Algebra {\bf 31} (2003), no.4, 1961-1973.


\bibitem{TS}
D. Takahashi, J. Satsuma {\it A soliton cellular automaton.} J. Phys. Soc. Japan {\bf 59} (1990), 3514-3519.
\bibitem{HKT}
G. Hatayama, A. Kuniba, T. Takagi {\it  Soliton cellular automata associated with crystal bases.}  Nuclear Phys. {\bf B 577} (2000), no. 3, 619--645. Ê
\bibitem{HHIKTT}
G. Hatayama, K. Hikami, R. Inoue, A. Kuniba, T. Takagi, T. Tokihiro. {\it The $A\sp {(1)}\sb M$ automata related to crystals of symmetric tensors. } J. Math. Phys. {\bf 42} (2001), no. 1, 274--308.
\bibitem{HKOT}
G. Hatayama, A. Kuniba, M. Okado, T.Takagi {\it Combinatorial $R$ matrices for a family of crystals: $B\sp {(1)}\sb n,D\sp {(1)}\sb n,A\sp {(2)}\sb {2n}$, and $D\sp {(2)}\sb {n+1}$ cases. } J. Algebra {\bf 247} (2002), no. 2, 577--615.
\bibitem{KNY}
K. Kajiwara, M. Noumi, Y. Yamada {\it Discrete dynamical systems with $W(A\sb {m-1}\sp {(1)}\times A\sb {n-1}\sp {(1)})$ symmetry. } Lett. Math. Phys. {\bf 60} (2002), no. 3, 211--219.
\bibitem{Kir}
A.N. Kirillov  {\it Introduction to tropical combinatorics. }
Physics and combinatorics, 2000 (Nagoya), 82--150,
World Sci. Publishing, River Edge, NJ, 2001.
\bibitem{Y}
Y.Yamada {\it A birational representation of Weyl group,
combinatorial $R$-matrix and discrete Toda equation.}  Physics and
combinatorics, 2000 (Nagoya), 305--319, World Sci. Publishing,
River Edge, NJ, 2001.
\bibitem{KOTY}
A.Kuniba, M. Okado, T. Takagi, Y. Yamada. {\it Tropical $R$ and
tau functions.}  Comm. Math. Phys. {\bf 245} (2004), no. 3,
491--517.
\bibitem{NouY}
M. Noumi, Y. Yamada {\it  Tropical Robinson-Schensted-Knuth correspondence and birational Weyl group actions.}  Representation theory of algebraic groups and quantum groups, 371--442, Adv. Stud. Pure Math., {\bf 40}, Math. Soc. Japan, Tokyo, 2004.




\bibitem{V1}
A.P. Veselov  {\it Yang-Baxter maps and integrable dynamics.}  Physics Letters A, {\bf 314} (2003), 214-221.
\bibitem{GV}
V.M. Goncharenko, A.P. Veselov {\it Matrix solitons and Yang-Baxter maps.}
In "New Trends in Integrability and Partial Solvability."
Edited by A.B. Shabat et al, Kluwer Academic Publishers, 2004, pp. 191-197.
\bibitem{SV}
Yu.B. Suris, A.P. Veselov {\it Lax pairs for Yang-Baxter maps.} J. Nonlin. Math. Phys. {\bf 10} (2003), suppl. 2, 223--230.
\bibitem{RV}
N. Reshetikhin and A. Veselov {\it Poisson Lie groups and Hamiltonian theory of the Yang-Baxter maps.}
math.QA/0512328.


\bibitem{Yang}
C.N. Yang {\it Some exact results for the many-body problem in one dimension with repulsive delta-function interaction.}
Phys. Rev. Lett. 19 (1967), 1312-1315.
\bibitem{B}
R. Baxter {\it Exactly solvable models in statistical mechanics.} Acad. Press, 1982.
\bibitem{TF}
L.A. Takhtajan, L.D. Faddeev {\it The quantum inverse problem method and the $XYZ$ Heisenberg model.}
Uspekhi Mat. Nauk 34 (1979), no. 5(209), 13--63.
\bibitem{Gaudin}
M. Gaudin {\it La fonction d'onde de Bethe.} Masson, 1983.
\bibitem{JM}
M. Jimbo, T. Miwa {\it Algebraic analysis of solvable lattice models.} AMS, 1995.

\bibitem{VS}
A.P. Veselov, A.B. Shabat {\it Dressing chain and spectral theory of the Schr\"odinger operator.}
Funct Anal Appl. 27:2 (1993), 1-21.
\bibitem{A}
V.E. Adler {\it Recutting of polygons.} Funct. Anal. Appl. 27:2 (1993), 79-80.
\bibitem{Nov}
S.P. Novikov {\it Periodic problem for KdV equation. I} Funct Anal
Appl. {\bf 8:3} (1974), 54-66.
\bibitem{DMN}
B.A. Dubrovin, V.B. Matveev, S.P.Novikov {\it KdV type equations
and finite-gap operators.} Russian Math Surveys (1976) no. 1(187), 55--136.




\bibitem{V}
A.P. Veselov {\it Integrable mappings.} Russian Math. Surveys, 46:5 (1991), 1-51.

\bibitem{Kul}
P.P. Kulish {\it Factorization of the classical and quantum S-matrix and conservation laws.}
Theor. Math. Phys. 26 (1976), 132-137.
\bibitem{Sol}
S. P. Novikov, S. V. Manakov,  L. V. Pitaevski, V. E. Zakharov,  {\it
Soliton theory: The Inverse Scattering Method.}  Plenum: New York (1984).
\bibitem{G}
V.M. Goncharenko {\it Multisoliton solutions of the matrix KdV equation.}
Theor. Math. Phys., 126:1 (2001), 81-91.
\bibitem{Tsu}
T. Tsuchida {\it $N$-soliton collision in the Manakov model.} Prog. Theor. Phys. {\bf 111} (2004), no.2, 151-182.
\bibitem{APT}
M. Ablowitz, B. Prinari, A.D. Trubatch {\it Soliton interactions in the vector NLS equation.}
Inverse Problems {\bf 20} (2004), no. 4, 1217--1237.

\bibitem{ABS2}
V.E. Adler, A.I. Bobenko, Yu. B.  Suris {\it Geometry of
Yang-Baxter maps: pencils of conics and quadrirational mappings. }
Comm. Anal. Geom. {\bf 12} (2004), no. 5, 967--1007.
\bibitem{ABS}
V.E. Adler, A.I. Bobenko, Yu.B. Suris. {\it Classification of
integrable equations on quad-graphs. The consistency approach.}
Comm. Math. Phys., 2003, {\bf 233}, 513-543.
\bibitem{BS}
A.I. Bobenko, Yu.B. Suris. {\it Integrable systems on
quad-graphs.} Int. Math. Res. Notices, 2002, No. 11, 573-611.
\bibitem{N}
F.W. Nijhoff. {\it Lax pair for the Adler (lattice
Krichever-Novikov) system.} Phys. Lett. A, 2002, {\bf 297}, 49-58.



\bibitem{FR}
I.B. Frenkel, N. Yu. Reshetikhin {\it Quantum affine algebras and holonomic difference equations.}
Comm. Math. Phys. 146 (1992), 1-60.
\bibitem{FK}
S. Fomin, A.N. Kirillov {The Yang-Baxter equation, symmetric functions and Schubert polynomials.}
Discrete Math. 153 (1996), 123-143.
\bibitem{BFZ}
A. Berenstein, S. Fomin, A. Zelevinsky {\it Parametrizations of canonical bases and totally positive matrices.}
Adv. in Math. 122 (1996), 49-149.
\bibitem{NY}
M. Noumi, Y. Yamada {\it Affine Weyl groups, discrete dynamical
systems and Painlev\'e equations.} Comm. Math. Phys. 199 (1998),
no. 2, 281--295.

\bibitem{Bour}
N. Bourbaki {\it Groupes et alg\`ebres de Lie.} Chap.~VI, Masson, 1981.

\bibitem{GPV}
M.-P. Grosset, V. Papageorgiou, A.P. Veselov {\it Periodic problem for discrete KdV equation and Yang-Baxter correspondences.} In preparation.

\bibitem{GLR}
I. Gohberg, P. Lancaster, L. Rodman {\it Matrix polynomials.} New York: Academic Press, 1982.



\bibitem{MV}
J. Moser, A.P. Veselov {\it Discrete versions of some classical integrable systems
and factorization of matrix polynomials.}
Comm. Math. Phys., 139 (1991), 217-243.
\bibitem{R}
N. Reshetikhin {\it Integrable discrete systems.}  Quantum groups and their applications in physics
(Varenna, 1994), 445-487, Proc. Internat. School Phys. Enrico Fermi, 127, IOS, Amsterdam, 1996.
\bibitem{HKKR}
T. Hoffmann, J. Kellendonk, N. Kutz, N. Reshetikhin {\it Factorization dynamics and Coxeter-Toda lattices.} Comm. Math. Phys. 212 (2000), 297-321.
\bibitem{RS}
N.Yu. Reshetikhin, M.A. Semenov-Tian-Shansky {\it Quantum R-matrices and factorisation problems.}
J. Geom. Phys. 5 (1988), 533-550.
\bibitem{O}
A.Odesskii {\it Set-theoretical solutions to the Yang-Baxter relation from factorization of matrix polynomials and theta-functions.}
Mosc. Math. J. {\bf 3} (2003), no. 1, 97--103.


\bibitem{Lax}
P.\,~D.\,~Lax. {\it Integrals of nonlinear equations of evolution
and solitary waves. } {\sl  Comm. Pure Appl. Math. ,} {\bf 21,}
467-490 (1968).

\bibitem{RSTS}
A.G. Reyman, M.A. Semenov-Tian-Shanski {\it Integrable Systems.
Group Theoretic Approach.}, Modern Mathematics, RCD, Izhevsk, 2003.
\bibitem{Drin}
V.G. Drinfeld {\it Hamiltonian structures on Lie groups, Lie
bialgebras and geometrical meaning of the Yang-Baxter equation.}
Sov. Math. Dokl., 268 (1983), 285-287.
\bibitem{W}
A. Weinstein {\it  Poisson geometry. Symplectic geometry.}
Differential Geom. Appl.  9  (1998),  no. 1-2, 213--238.


\bibitem{Berger} 
M.~Berger {\it Geometry.} Springer-Verlag, Berlin, 1987.
\bibitem{AS}
V.E. Adler, Yu. B. Suris ${\rm Q}\sb 4$: {\it integrable master equation related to an elliptic curve.}  Int. Math. Res. Not. 2004, no. 47, 2523--2553.
\bibitem{PTV}
V. Papageorgiou, A. Tongas and A.P. Veselov {\it Yang-Baxter maps and symmetries of integrable equations on quad-graphs.} math.QA/0605206. To appear in J. Math. Phys. (2006).





\end{thebibliography}
\end{document}